\documentclass[twoside,11pt]{article}
\usepackage{graphicx}
\usepackage{subfig}
\usepackage{amsmath}
\usepackage{amssymb}
\usepackage{amsthm}
\usepackage{url}
\newcommand\nonu{\nonumber}
\newcommand\sa{\smallskipamount}
\newcommand\ma{\medskipamount}
\newcommand\ba{\bigskipamount}
\newcommand\mLP{\\[\ma]}
\newcommand\bPP{\\[\ba]\indent}
\newcommand\bLP{\\[\ba]}
\newcommand\sPP{\\[\sa]\indent}
\newcommand\CC{\mathbb{C}}
\newcommand\RR{\mathbb{R}}
\newcommand\ZZ{\mathbb{Z}}
\newcommand\ep\varepsilon
\newcommand\Ga{\Gamma}
\newcommand\half{\frac12}
\newcommand\id{\operatorname{id}}
\newcommand\iy\infty
\newcommand\const{{\rm const.}\,}
\newcommand{\hyp}[5]{\,\mbox{}_{#1}F_{#2}\left(
  \genfrac{}{}{0pt}{}{#3}{#4};#5\right)}
\newcommand\RHS{right-hand side}
\renewcommand\Re{\operatorname{Re}}
\newcommand\ri[1]{\genfrac{}{}{0pt}{}{#1}{\longrightarrow}}
\newcommand\li[1]{\genfrac{}{}{0pt}{}{#1}{\longleftarrow}}
\numberwithin{equation}{section}
\newtheorem{theorem}{Theorem}[section]
\newtheorem{corollary}[theorem]{Corollary}
\newtheorem{Remark}[theorem]{Remark}
\newenvironment{remark}{\begin{Remark}\rm}{\end{Remark}}
\newcommand\Proof{\noindent{\bf Proof}\quad}

\usepackage{fancyhdr}
\begin{document}
\title{Integral representations for Horn's $H_2$ function\\
and Olsson's  $F_P$ function}
\author{Enno Diekema and Tom H. Koornwinder}
\date{}
\maketitle
\begin{abstract}
We derive some Euler type double integral representations for
hypergeometric functions in two variables. In the first part of this
paper we deal with Horn's $H_2$ function, in the second part with
Olsson's $F_P$ function. Our double integral representing the $F_P$
function is compared with the formula for the same integral
representing an $H_2$ function by M. Yoshida (Hiroshima Math.\ J.\
10 (1980), 329--335) and M. Kita (Japan.\ J.\ Math.\ 18 (1992), 25--74). 
As specified by Kita, their integral is defined by a
homological approach.  We present a classical double integral version
of Kita's integral, with outer integral over a Pochhammer double loop,
which we can evaluate as $H_2$ just as Kita did for his integral.
Then we show that shrinking of the double loop yields a sum of two
double integrals for $F_P$.
\end{abstract}
\bLP
\section{Introduction}
This paper deals with integral representations for two bivariate
hypergeometric functions: on the one hand
the $H_2$ Horn function \cite[(5.7.14)]{6} which occurs in Horn's
list \cite{9} of bivariate hypergeometric series of order two;
on the other hand the $F_P$ function which was
introduced by Olsson \cite{12}, \cite{13}. The integral representations
under consideration will be integrals of classical type over
two-dimensional real domains.
Necessarily,
because of convergence, this requires constraints on the allowed
parameter values.

Horn \cite{8} called a
double power series
$\sum_{i,j=0}^\iy A(i,j) x^i y^j$
{\em hypergeometric} if the two quotients $\frac{A(i+1,j)}{A(i,j)}$
and $\frac{A(i,j+1)}{A(i,j)}$ are rational functions of $i$ and $j$.
Write the two rational functions as quotients of polynomials in $i$
and $j$ without common factors:
\[
\frac{A(i+1,j)}{A(i,j)}=\frac{F(i,j)}{F'(i,j)}\qquad
{\rm and}
\qquad\frac{A(i,j+1)}{A(i,j)}=\frac{G(i,j)}{G'(i,j)}\,.
\]
In addition it is assumed that
$F'(i,j)$ contains the factor $i+1$ and
$G'(i,j)$ the factor $j+1$.
Then the highest degree in $i,j$ of the four polynomials $F,F',G,G'$
is called the {\em order} of the
hypergeometric series.

Horn \cite{9} listed all convergent bivariate hypergeometric series
of order two. For the case that all polynomials 
$F,F',G,G'$ have degree two he obtains 14 items, which include
Appell's hypergeometric series $F_1$, $F_2$, $F_3$ and $F_4$,
dating back to 1880 \cite{2}, \cite{19}, \cite{3}, but also,
among others,
the $H_2$ function. These bivariate hypergeometric series can be
considered as the natural two-variable analogues of the Gauss
hypergeometric series. Horn gives 20 further items which are
confluent cases of the first 14 items.
For the basics of hypergeometric functions of two variables
see for instance \cite[Ch.~9]{4}, \cite[Sections 5.7--5.12]{6},
\cite[Ch.~1]{7} and \cite[Ch.~8]{15}.

In \cite[Section 5.9]{6},
for each function in Horn's list a system of two partial differential
equations (PDEs) is given
which has this function as a solution. This follows a method already
indicated by Horn in \cite {8} and \cite{9}. Some of these systems
can be transformed into each other. This is
for instance the case with the systems for $F_2$, $F_3$ and $H_2$.
Otherwise stated, $F_3$ and $H_2$, when suitably adapted, occur
as solutions of the system \cite[5.9(10)]{6} 
for $F_2$. A comprehensive list of solutions of this $F_2$ system
was given by Olsson \cite[p.~1289, Table I]{14}. Each solution is in
terms of the six bivariate functions $F_2$, $F_3$, $H_2$, $F_P$,
$F_Q$ and $F_R$. While $F_R$ is not even a hypergeometric series
in the sense of Horn, $F_P$ and $F_Q$ are hypergeometric series
of order three. 

Euler type double integral representations for the Appell
hypergeometric functions $F_1,F_2,F_3$
were already given by Appell \cite[Chap.~III]{19} in 1882.
Some such representations
for other functions in Horn's list are scattered through the literature.
Kita \cite[pp.~56--58]{10} gives a long list of double integral
representations, however with integrals
defined by a homological approach.
Single integral representations with one or two hypergeometric
functions in the integrand can be found in the literature in more
cases. Usually double integral representations imply such single
integral representations.

When we concentrate on solutions of the $F_2$ system, it is natural to ask if, beside $F_2$
and $F_3$, also other solutions have an Euler type double integral
representation. In this paper we obtain these for $H_2$ and $F_P$.
The integral representation we derive for $H_2$ turned out to
have been given earlier in a forgotten paper by Tuan \& Kalla
\cite[(77)]{16}. We show that our method of derivation gives rise to
several variants of this formula.

The integral representation we derive for $F_P$ is puzzling when
compared with Yoshida \cite[(0.10)]{17} and
Kita \cite[p.57, item 9]{10}. A straightforward rewriting of our
integral as an integral over a triangular region is evaluated
by these authors as an $H_2$ function. However,
as becomes clear from Kita's paper, their double integral is
defined by a homological approach, involving twisted cycles.
While Kita gives, for instance,
an integral of this type for $F_2$ which can be
immediately matched with Appell's classical integral representation
\cite[5.8(2)]{6} (under suitable
constraints on the parameter values), a similar relation apparently
does not exist in the $F_P$ case between Kita's and our formula.
In this paper we will not go into the technicalities of the
homological approach (see for instance
\cite[\S3]{10}, \cite[\S2.3]{20}, \cite[Ch.~IV]{23}).
Instead, we will show
that the evaluation as $H_2$ by Yoshida and Kita
can also be achieved by a classical double integral, where the outer
integral is taken over a Pochhammer double loop. Then it is clear
that a branch point caused by an independent variable prevents
the Pochhammer loop from being shrunk such that there results
an integral over a segment.
However, it is possible to do the shrinking in such a way
that one arrives at a sum of two integrals
over real intervals. Both integrals in this sum can be evaluated as
an $F_P$. The resulting formula is a three-term relation
\cite[(53)]{14} involving one $H_2$ and two $F_P$ terms.
A similar phenomenon, but more simple, can already be illustrated
for a variant of the
Euler integral representation for the Gauss hypergeometric function.
While it is quite standard to derive connection formulas by
using the homological approach (see Remark \ref{th5}),
this does not seem to have been done before by using integrals
over Pochhammer double loops.

The contents of this paper are as follows.
In Section \ref{s2} we summarize the integral representations of the
Gauss hypergeometric function. In Section \ref{s3} we derive some
integral representations of the $H_2$ function. We also indicate how
to get many further integral representations of this function. In
Section \ref{s4} we derive a double integral representation of the
$F_P$ function and we give some corollaries.
In Section \ref{s5} we give
another corollary of this double integral representation for~$F_P$.
We compare it with the
Yoshida-Kita integral representation which uses the homological
approach and yields~$H_2$. 
Then we show that this $H_2$ evaluation
already occurs with a classical double integral, the outer one
being over a Pochhammer double loop. Shrinking of the double loop 
leads to a known three-term relation.
Finally, in Section \ref{s6} we briefly discuss how the obtained double
integrals appear when written as integral representations of solutions
of the system of PDEs for~$F_2$.
\paragraph{Notation}
For $a\in\CC$ and $k\in\ZZ$ the {\em Pochhammer symbol} is defined by
\begin{equation*}
(a)_k:=\frac{\Ga(a+k)}{\Ga(a)}=
\begin{cases}
a(a+1)\ldots(a+k-1),&k>0,\\
1,&k=0,\\
(-1)^k/(1-a)_{-k}\,,&k<0.
\end{cases}
\end{equation*}

\section{Euler integral representations for the Gauss
hypergeometric function}
\label{s2}
The {\em Gauss hypergeometric function} \cite[Ch. 2]{1} is defined
as a power series by
\begin{equation}
\hyp21{a,b}{c}z:=
\sum_{k=0}^\iy\frac{(a)_k(b)_k}{(c)_k\,k!}\,z^k
\qquad(|z|<1).
\label{2.1}
\end{equation}
It has, as a function of $z$, a one-valued analytic continuation to
$\CC\backslash[1,\iy)$, as is seen (at least under the given parameter
restrictions) from the {\em Euler integral representation}
\begin{align}
\hyp21{a,b}{c}z
&=\frac{\Ga(c)}{\Ga(b)\Ga(c-b)}\,
\int_0^1 t^{b-1}(1-t)^{c-b-1}(1-zt)^{-a}\,dt\label{2.2}\\
&=\frac{\Ga(c)}{\Ga(b)\Ga(c-b)}\,(1-z)^{c-a-b}
\int_0^1 t^{c-b-1}(1-t)^{b-1}(1-zt)^{a-c}\,dt\label{2.3}\\
&=\frac{\Ga (c)}{\Ga(a)\Ga(c-a)}\,
\int_0^1 t^{a-1}(1-t)^{c-a-1}(1-zt)^{-b}\,dt\label{2.4}\\
&=\frac{\Ga(c)}{\Ga(a)\Ga(c-a)}\,(1-z)^{c-a-b}
\int_0^1 t^{c-a-1}(1-t)^{a-1}(1-zt)^{b-c}\,dt\label{2.5}
\end{align}
with convergence conditions $\Re c>\Re b>0$
for \eqref{2.2} and \eqref{2.3},
and $\Re c>\Re a>0$ for \eqref{2.4} and \eqref{2.5}.
Formula \eqref{2.2} is the usual formula for the Euler integral representation.
Then \eqref{2.3} follows from \eqref{2.2} by the
{\em Euler transformation formula}
\begin{equation}
\hyp21{a,b}{c}z=(1-z)^{c-a-b}\hyp21{c-a,c-b}{c}z.
\label{2.6}
\end{equation}
Note its special case
\begin{equation}
\hyp21{a,b}{a}z=(1-z)^{-b}.
\label{eq33}
\end{equation}
Formulas \eqref{2.4} and \eqref{2.5} follow by the symmetry of
\eqref{2.1} in $a$ and $b$.
The {\em Pfaff transformation formula} 
\begin{equation}
\hyp21{a,b}{c}z=(1-z)^{-a} \hyp21{a,c-b}c{\frac z{z-1}}
\label{2.7}
\end{equation}
follows from \eqref{2.2} or \eqref{2.5} by the change of
integration variable
$t\rightarrow 1-t$. Similarly, its variant 
\begin{equation*}
\hyp21{a,b}{c}z=(1-z)^{-b} \hyp21{c-a,b}c{\frac z{z-1}}
\label{2.8}
\end{equation*}
follows from \eqref{2.3} or \eqref{2.4}. For generic values of the
parameters, formulas \eqref{2.2}--\eqref{2.5} are the only integral
representations of the Gauss hypergeometric function of the form
\[
\const(1-z)^\lambda\int_{0}^{1}t^\alpha(1-t)^\beta(1-zt)^\Ga\,dt,
\]
or else there would be other transformation formulas of the form
\eqref{2.6}.

On the \RHS\ of \eqref{2.2} we can apply a one-parameter group of
transformations of integration variable
\begin{equation}
\phi_p\colon t\longmapsto \frac{t}{p+(1-p)t}\qquad(p>0),
\label{2.9}
\end{equation}
which map the integration interval $[0,1]$ onto itself.
Note that $\phi_p\circ \phi_q=\phi_{pq}$ and $\phi_1=\id$. The
integral representation resulting from applying the transformation of
integration variable $\phi_p$ to \eqref{2.2} is
\begin{align}
\hyp21{a,b}cz
&=\frac{\Ga(c)}{\Ga(b)\Ga(c-b)}\,p^{-b}\nonu\\
&\qquad\times\int_0^1 t^{b-1}(1-t)^{c-b-1}
\left(1-\frac{p-1}{p}t\right)^{a-c}
\left(1-\frac{z+p-1}{p}t\right)^{-a}dt\label{2.10}\\
&=p^{-b}\,F_1\left(b,c-a,a,c;\frac{p-1}{p},\frac{z+p-1}{p}\right),
\label{2.11}
\end{align}
where the constraints in \eqref{2.10} are the same as in \eqref{2.2},
and where the second equality follows from the integral representation
\cite[5.8(5)]{6} for the Appell hypergeometric function $F_1$.  Then
\eqref{2.11} gives a reduction formula for $F_1( a,b,c,b+c;x,y)$,
which is \cite[5.10(1)]{6} combined with~\eqref{2.7}.

Similarly, application of the transformation of integration variable
\begin{equation}
\psi_p(t):=\phi_p(1-t)=\frac{1-t}{1-(1-p)t}\qquad(p>0)
\label{2.12}
\end{equation}
to the \RHS\ of \eqref{2.2} gives 
\begin{align}
\hyp21{a,b}cz
&=\frac{\Ga(c)}{\Ga(b)\Ga(c-b)}\,p^{c-b}(1-z)^{-a}\nonu\\
&\qquad\times\int_0^1 t^{c-b-1}(1-t)^{b-1}\big(1-(1-p)t\big)^{a-c}
\left( 1-\frac{1-p-z}{1-z}t\right)^{-a}dt\label{2.13}\\
&=p^{c-b}(1-z)^{-a}\,
F_1\left(c-b,c-a,a,c;1-p,\frac{1-p-z}{1-z}\right),
\label{2.14}
\end{align}
again with the same constraints in \eqref{2.13} as in \eqref{2.2}, and
with \eqref{2.14} being the reduction formula \cite[5.10(1)]{6}
combined with \eqref{2.6}. Note that the special case $p=1$ of
\eqref{2.13} or \eqref{2.14} gives \eqref{2.7}.

\begin{remark}
\label{rem2.1}
The same method of applying $\phi_p$ or $\psi_p$ to the
integration variable works for the single integral representation
\cite[5.8(5)]{6} for the Appell $F_1$ function or, more generally,
for the single integral representation 
\cite[p.116, formula (8)]{3}
for the Lauricella $F_D$
function in $n$ variables.
For instance, in the Lauricella case the result
is an integral representation for Lauricella $F_D$ in
$n+1$ variables, by which we also obtain a reduction formula
for Lauricella $F_D$ with a dependence between the parameters.
\end{remark}
\section{Integral representations for the $H_2$ function}
\label{s3}
\begin{figure}[ht]
\centering
\parbox{5cm}{
\includegraphics[width=5cm]{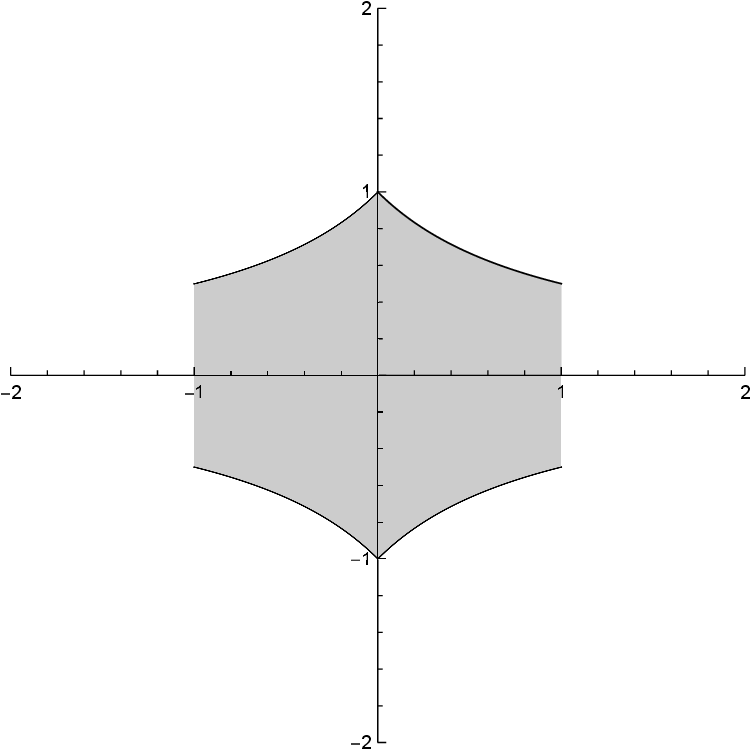}
\caption{Convergence region $\Omega_1\cap\RR^2$
of the $H_2$ function.}
\label{fig:1}}
\qquad
\begin{minipage}{5cm}
\includegraphics[width=5cm]{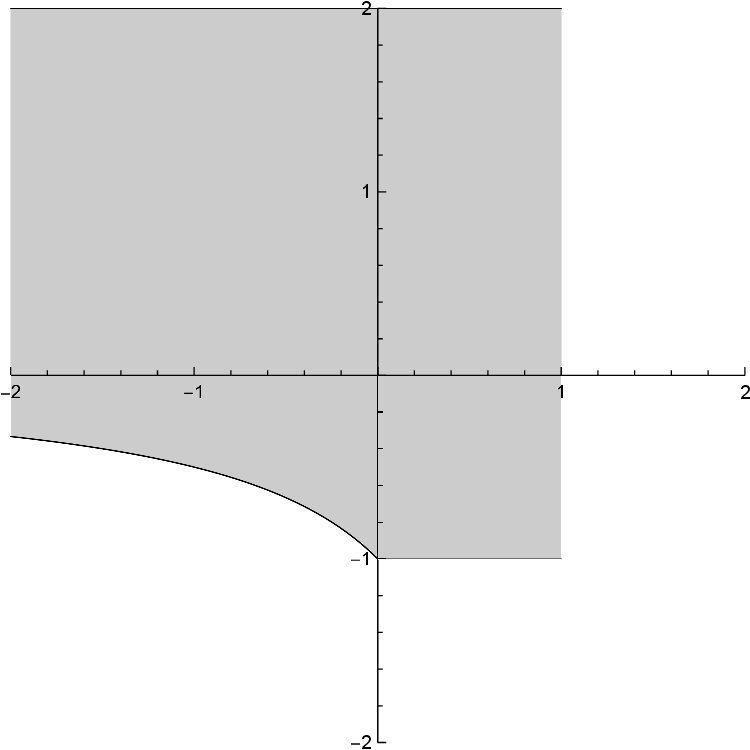}
\caption{Region of analytic continuation $\Omega_2\cap\RR^2$
of $H_2$.}
\label{fig:2}
\end{minipage}
\end{figure}

We start with the double power series of the $H_2$ function
\cite[5.7(14) and (49)]{6}
\begin{equation}
H_2(a,b,c,d,e;x,y)=H_2(x,y):=
\sum_{i,j=0}^\iy
\frac{(a)_{i-j}(b)_{i}( c)_{j}(d)_{j}}{(e)_{i}\,i!\,j!}\,x^i y^j
\label{3.1}
\end{equation}
with convergence region (see Figure \ref{fig:1})
\[
\Omega_1:=\big\{(x,y)\in\CC^2\;\big|\; |x|<1,\;|y|<(|x|+1)^{-1}\big\}.
\]
By inserting 
\[
(a)_{i-j}=(-1)^{j}\,\frac{(a-j)_{i}}{(1-a)_{j}}
\]
in \eqref{3.1} we get
\begin{align}
H_2(a,b,c,d,e;x,y)
&=\sum_{j=0}^\iy\frac{(c)_{j}(d)_{j}}{(1-a)_{j}\,j!}\,(-y)^{j}
\sum_{i=0}^\iy\frac{(a-j)_{i}(b)_{i}}{(e)_{i}\,i!}\,x^{i}
\nonu\\
&=\sum_{j=0}^\iy\frac{(c)_{j}(d)_{j}}{(1-a)_{j}\,j!}\,(-y)^j
\hyp21{a-j,b}{e}x
\label{3.2}
\end{align}
for $(x,y)\in\Omega_1$. Specialization of \eqref{3.2} gives
\begin{equation}
H_2(a,b,c,d,e;0,y)=\hyp21{c,d}{1-a}{-y},
\label{eq35}
\end{equation}
initially for $|y|<1$, and after analytic continuation for complex $y$ outside
$(-\iy,-1]$.

Substitution of Euler's integral representation \eqref{2.2}
in \eqref{3.2} gives
\[
H_2(a,b,c,d,e;x,y)
=\sum_{j=0}^\iy\frac{(c)_{j}(d)_{j}}{(1-a)_{j}\,j!}\,(-y)^{j}\,
\frac{\Ga(e)}{\Ga(b)\Ga(e-b)}
\int_{0}^{1}u^{b-1}(1-u)^{e-b-1}(1-xu)^{-a+j}\,du
\]
with condition $0<\Re b<\Re e$. Because,
for $(x,y)\in\Omega_1$, we have
\begin{multline*}
\sum_{j=0}^\iy\,
\left|\frac{(c)_{j}(d)_{j}}{(1-a)_{j}\,j!}\,(-y)^{j}\right|
\int_{0}^{1}\left|u^{b-1}(1-u)^{e-b-1}( 1-xu)^{-a+j}\right|\,du\\
\le\const\sum_{j=0}^\iy\,
\left|\frac{(c)_{j}(d)_{j}}{(1-a)_{j}\,j!}\right|\;
|y|^j \big(|x|+1\big)^{j}
\int_{0}^{1}\left| u^{b-1}(1-u)^{e-b-1}(1-xu)^{-a}\right|\,du<\iy\,,
\end{multline*}
the sum and integral can be interchanged by dominated convergence. So
\begin{align}
&H_2(a,b,c,d,e;x,y)\nonu\\
&\qquad=\frac{\Ga(e)}{\Ga(b)\Ga(e-b)}
\int_{0}^{1}u^{b-1}(1-u)^{e-b-1}(1-xu)^{-a}
\sum_{j=0}^\iy\frac{(c)_{j}(d)_{j}}{(1-a)_{j}\,j!}\,
(-y)^{j}(1-xu)^{j}\,du
\nonu\\
&\qquad
=\frac{\Ga(e)}{\Ga(b)\Ga(e-b)}
\int_{0}^{1} u^{b-1}(1-u)^{e-b-1}(1-xu)^{-a}
\hyp21{d,c}{1-a}{-y(1-xu)}du.
\label{3.3}
\end{align}
\newline
The absolute value of the integrand of \eqref{3.3} is dominated by
$\const\big|u^{b-1}(1-u)^{e-b-1}\big|$, uniformly for $(x,y)$ in
compact subsets of the region
\[
\Omega_2: =\left\{ \left( x,y\right) \in \CC^{2}\,\big|\,
x\notin[1,\iy)\,\mbox{ and }\;\forall u\in[0,1]\;y(xu-1) \notin[1,\iy)\right\}.
\]
Therefore, under the constraints $0<\Re b<\Re e$, \eqref{3.3}
provides an analytic continuation of $H_2(x,y)$
for $(x,y)\in \Omega_2$.
The intersection of $\Omega_2$ with $\RR^2$  is
\[
\left\{(x,y)\in\RR^2\;\big|\;x\le 0,\;( x-1) y<1\;\;\vee\;\;
0<x<1,\;y>-1 \right\}
\]
(see Figure \ref{fig:2}), much larger than the convergence
region for the power series \eqref{3.1} intersected with $\RR^2$
(see Figure \ref{fig:1}).

Substitution of Euler's integral representation \eqref{2.2} in
\eqref{3.3} gives
\begin{multline}
H_2(a,b,c,d,e;x,y)
=\frac{\Ga(e)}{\Ga(b)\Ga(e-b)}\,\frac{\Ga(1-a)}{\Ga(c)\Ga(1-a-c)}\\
\times\int_{0}^{1}u^{b-1}(1-u)^{e-b-1}(1-xu)^{-a}
\int_{0}^{1}v^{c-1}(1-v)^{-a-c}(1+yv-xyuv)^{-d}\,dv\,du,
\label{3.4}
\end{multline}
where we require moreover that $0<\Re c<\Re(1-a)$.
By Fubini's theorem formula \eqref{3.4} can be rewritten as
\begin{multline}
H_2(a,b,c,d,e;x,y)=\frac{\Ga(e)}{\Ga(b)\Ga(e-b)}\,
\frac{\Ga(1-a)}{\Ga(c)\Ga(1-a-c)}\\
\times\int_{0}^{1}\int_{0}^{1}u^{b-1}v^{c-1}(1-u)^{e-b-1}
(1-v)^{-a-c}(1-xu)^{-a} (1+yv-xyuv)^{-d}\,du\,dv,
\label{3.5}
\end{multline}
because the absolute value of the integrand is dominated by\quad
$\const |u^{b-1}(1-u)^{e-b-1}v^{c-1}\times$ $(1-v)^{-a-c}|$, for which
the integral over $(u,v) \in[ 0,1]\times [0,1]$ is finite. So formula
\eqref{3.5} is valid for $x\in \Omega_2$ and $0<\Re b<\Re e$,
$0<\Re c<\Re(1-a)$.

Formula (3.3) was earlier given by Olsson \cite[(44)]{14}. The
integral representation \eqref{3.5} was earlier given by Tuan \& Kalla
\cite[(77)]{16}. In both references the same constraints on the
parameters are given as we have, but the region on which the integral
representation is valid is not specified in these two references. For
the proof of \cite[(77)]{16} Tuan \& Kalla just mention that the
formula follows by termwise integration of the double power series of
the integrand. This is essentially the proof as given above.

For $H_2$ we have the transformation formula
\begin{equation}
H_2(a,b,c,d,e;x,y)=(1-x)^{-a}H_2\left(a,e-b,c,d,e;\frac{x}{x-1},y(1-x)\right),
\label{3.6}
\end{equation}
which is a rewriting of \cite[(25)]{14}. Formula \eqref{3.6}
also follows
from \eqref{3.2} together with \eqref{2.7} or from \eqref{3.3} or
\eqref{3.5} by the change of integration variable $u\rightarrow 1-u$.

In the proof of \eqref{3.5} we have twice substituted Euler's integral
representation \eqref{2.2} for a $_{2}F_{1}$: after \eqref{3.2} and
after \eqref{3.3}. In both places we might have substituted one of the
alternatives \eqref{2.3}, \eqref{2.4} or \eqref{2.5} for
\eqref{2.2}. After \eqref{3.2} only the alternative \eqref{2.3} will
be an option, because $\Re(a-j)$ would not be positive for $j$ large
enough and therefore the conditions for \eqref{2.4} and \eqref{2.5}
would be violated.

After substitution in \eqref{3.2} of \eqref{2.3} we can proceed in the
same way as above for the derivation of (2.5). We obtain:
\begin{align}
&H_2(a,b,c,d,e;x,y)=\frac{\Ga(e)}{\Ga(b)\Ga(e-b)}\,
(1-x)^{e-a-b}\nonu\\
&\qquad\qquad\qquad\qquad\times
\int_{0}^{1} u^{e-b-1}(1-u)^{b-1}(1-xu)^{a-e}
\hyp21{d,c}{1-a}{\frac{y(x-1)}{1-xu}} du
\label{3.7}\\
\noalign{\allowbreak}
&\qquad=\frac{\Ga(e)}{\Ga(b)\Ga(e-b)}\,\frac{\Ga(1-a)}{\Ga(c)\Ga(1-a-c)}
\int_{0}^{1}\int_{0}^{1}u^{e-b-1}(1-u)^{b-1}v^{c-1}(1-v)^{-a-c}\nonu\\
&\qquad\qquad\qquad\qquad\qquad\qquad\qquad\qquad
\times(1-xu)^{a+d-e}\big(1-xu+y(1-x)v\big)^{-d}\,du\,dv.
\label{3.8}
\end{align}
Formulas (3.7) and (3.8) give an analytic continuation of
$H_2(x,y)$ to the region
\[
\Omega_3:=\left\{(x,y) \in\CC^2\,\big|\,
x\notin [1,\iy)\;\mbox{ and }\;
\forall u\in[0,1]\;\frac{y(x-1)}{1-xu}\notin [1,\iy)\right\}
\]
(note that $\Omega_3$ has the same intersection with 
$\RR^2$ as $\Omega_2$). The parameter constraints for \eqref{3.7} and 
\eqref{3.8} are the same as those for \eqref{3.3} and \eqref{3.5}, respectively.

We can obtain two other variants of the double integral
representations \eqref{3.5} and \eqref{3.8} by substituting
\eqref{2.3} in \eqref{3.3} and \eqref{3.7}, respectively, but we omit
the explicit formulas. On the other hand, nothing of interest comes
out when we substitute \eqref{2.4} or \eqref{3.5} in \eqref{3.3} or
\eqref{3.7}. The resulting formulas only differ from the four earlier
integral representations by interchange of $c$ and $d$. This is the
trivial symmetry which is obvious from (3.1).

We may also vary the four double integral representations for $H_2$ by
applying the change of integration variable $\phi_{p}$ or $\psi_{p}$
(see \eqref{2.9}, \eqref{2.12}) to $u$ or $v$. As already remarked,
$\psi_{1}$ applied to $u$ corresponds to the transformation formula
\eqref{3.6}. Other cases possibly give specializations of double
integral representations for hypergeometric functions of three
variables (see also Remark~\ref{rem2.1}). Indeed, $\phi_p$ or $\psi_p$
applied to the $u$ variable in \eqref{3.5} would yield a double
integral with a relation between the parameters which corresponds to
the result of applying \eqref{2.11} or \eqref{2.14} to the ${}_2F_1$
in \eqref{3.2}.  If we were next to substitute for the resulting $F_1$
in \eqref{3.2} its double power series expansion then a hypergeometric
triple series with a relation between the parameters would arise which
is equal to~$H_2$.
The identities below, where we have omitted the constraints,
give the explicit result for the case that $\phi_p$ is applied.
\begin{align*}
&H_2(a,b,c,d,e;x,y)=\frac{\Ga(e)}{\Ga(b)\Ga(e-b)}\,
\frac{\Ga(1-a)}{\Ga(c)\Ga(1-a-c)}\,p^{-b}
\int_0^1\int_0^1 u^{b-1} v^{c-1} (1-u)^{e-b-1} (1-v)^{-a-c}
\\
&\qquad\qquad\qquad\qquad\qquad\qquad\qquad\qquad
\times\big(1-(1-p^{-1})u\big)^{a+d-e}
\big(1-p^{-1}(x+p-1)u\big)^{-a}\\
&\qquad\qquad\qquad\qquad\qquad\qquad\qquad\qquad
\times\big(1-(1-p^{-1})u+yv-p^{-1}(x+p-1)yuv\big)^{-d} du\,dv
\\
&\qquad\qquad\qquad\qquad
=p^{-b}\sum_{i,j,k=0}^\iy
\frac{(a)_{k-i}(b)_{j+k}(e-a)_{i+j}(c)_i( d)_i}
{(e)_{j+k}(e-a)_i\,i!\,j!\,k!}\,
y^i\left(\frac{p-1}p\right)^j\left(\frac{x+p-1}p\right)^k.
\end{align*}

Quite possibly the second identity can be generalized to a case
with less dependence between the parameters in the double integral
and in the triple series.
\section{Integral representations for the $F_P$ function}
\label{s4}
The function $F_P$ was first introduced by Olsson \cite{12}, \cite{13}
as a certain solution of the system of PDEs associated with Appell's
hypergeometric function $F_2$ given in \cite[5.9(10)]{6}.
In \cite{12} the notation $Z_1$ instead of $F_P$ is used.
The definitions \cite[(4)]{12} and \cite[(1a)]{13} of $F_P$ are also
different. But we can put together formulas for $F_P$ in \cite{12} and
\cite{13} as the following string of equalities:
\begin{align}
&F_P(a,b_1,b_2,c_1,c_2;x,y)\nonu\\
&\qquad=
\sum_{i,j=0}^\iy
\frac{(a)_{i+j}(a-c_{2}+1)_{i}(b_{1})_{i}(b_{2})_{j}}
{(a+b_{2}-c_{2}+1)_{i+j}(c_{1})_{i}\,i!\,j!}\,x^{i}(1-y)^j
\label{4.1}\\
\noalign{\allowbreak}
&\qquad=\sum_{i=0}^\iy
\frac{(a)_{i}(a-c_{2}+1)_{i}(b_{1})_{i}}
{(c_{1})_{i}(a+b_{2}-c_{2}+1)_{i}\,i!}\,
x^{i}\hyp21{a+i,b_{2}}{a+b_{2}-c_{2}+1+i}{1-y}
\label{4.2}\\
\noalign{\allowbreak}
&\qquad =y^{-a}\sum_{i=0}^\iy
\frac{(a)_{i}(a-c_{2}+1)_{i}(b_{1})_{i}}
{(c_{1})_{i}(a+b_{2}-c_{2}+1)_{i}\,i!}
\left(\frac{x}{y}\right)^i
\hyp21{a+i,a-c_{2}+1+i}{a+b_{2}-c_{2}+1+i}{\frac{y-1}{y}}
\label{4.3}\\
\noalign{\allowbreak}
&\qquad=y^{-a}\sum_{i,j=0}^\iy
\frac{(a)_{i+j}(a-c_{2}+1)_{i+j}(b_{1})_{i}}
{(a+b_{2}-c_{2}+1)_{i+j}(c_{1})_{i}\,i!\,j!}
\left(\frac{x}{y}\right)^i\left( \frac{y-1}{y}\right)^{j}.
\label{4.4}
\end{align}
The double power series \eqref{4.1} and \eqref{4.4} are both
of order three. By Horn's rule for determination of the
convergence region (see \cite[Section 5.7.2]{6})
the double power series \eqref{4.1} is absolutely convergent
on the region $|x|,|y-1|<1$ in $\CC^{2}$,
while \eqref{4.4} is absolutely convergent on the region
$|xy^{-1}|+|1-y^{-1}|<1$ in $\CC^{2}$,
which has intersection $|x|<1$, $y>\tfrac{1}{2}(1+|x|)$ with
$\RR^{2}$. Both convergence regions are neighborhoods of $(0,1)$. The
equalities \eqref{4.1}\,=\,\eqref{4.2} and \eqref{4.3}\,=\,\eqref{4.4}
on the respective convergence regions of \eqref{4.1} and \eqref{4.4}
follow by rewriting the double series with an inner $j$-sum and an
outer $i$-sum. The equality
\eqref{4.2}\,=\,\eqref{4.3}, which follows by Pfaff's transformation
\eqref{2.7}, is initially valid on the intersection of the two
convergence regions, but next yields an analytic continuation of $F_P$
to the union of the two convergence regions. The intersection of this
union with $\RR^{2}$ is $|x|<1\,\wedge\,y>0$.

By using the double power series \eqref{4.1} it can be seen that
$F_P(a,b_1,b_2,c_1,c_2;x,y)$ is the solution
of the system of PDEs for $F_2(a;b_1,b_2;c_1,c_2;x,y)$
which is regular and equal to $1$ at the point $(0,1)$.

In \cite{13} and \cite{14} Olsson takes \eqref{4.4}
as a definition of $F_P$, while he starts in \cite{12} with
\begin{equation}
F_P(a;b_{1},b_{2};c_{1},c_{2};x,y)=
\sum_{i=0}^\iy\frac{(a)_{i}(b_{2})_{i}}
{(a+b_{2}-c_{2}+1)_{i}\,i!}\,(1-y)^{i}
\hyp32{a+i,b_{1},a-c_{2}+1}{c_{1},a+b_{2}-c_{2}+1+i}x
\label{4.5}
\end{equation}
and then gives \eqref{4.1} and \eqref{4.2}.
Note that \eqref{4.5} follows on the convergence region of \eqref{4.1}
by rewriting \eqref{4.1} with an inner $j$-sum and an outer $i$-sum.

The expression of $F_P$ by the double power series \eqref{4.4}
implies the symmetry for $F_P$ which is given by the first equality
in \cite[(26)]{14}. However, no immediate symmetry can be derived from
\eqref{4.1} because the five parameters occur there in six shifted
factorials.

We now state our main result on $F_P$: a double integral representation.
\begin{theorem}
\label{th4.1}
For $\Re a,\Re(a-b_2),\Re(b_2+c_1-a),\Re(b_{2}-c_{2}+1),
\Re(c_{1}+c_{2}-a-1)>0$
and $(x,y)\in\CC^2$ such that
$x\notin[1,\iy)$ and $y\notin (-\iy,0]$, we have
\begin{multline}
F_P(a,b_{1},b_{2},c_{1},c_{2};x,y)=
\frac{\Ga(a+b_{2}-c_{2}+1)\Ga(c_{1})}
{\Ga(a)\Ga (a-c_{2}+1)\Ga(b_{2}-c_{2}+1)\Ga(c_{1}+c_{2}-a-1)}\\
\times y^{-b_2}\int_0^1\int_0^1 u^{a-1}
(1-u)^{b_{2}+c_{1}-a-1}v^{b_{2}-c_{2}}(1-v)^{c_{1}+c_{2}-a-2}
(1-xu)^{-b_{1}}\\
\times\big(u+(1-u)vy^{-1}\big)^{-b_2}\,du\,dv
\label{4.6}
\end{multline}
with absolutely convergent double integral.
\end{theorem}

\begin{remark}
In \eqref{4.6} the function of $(x,y)$ defined by the \RHS\ is
analytic on the given region and provides an analytic continuation of
the left-hand side. A similar remark applies to the subsequent
corollaries.
\end{remark}

Before proving Theorem \ref{th4.1}
we mention three immediate consequences, which are, in a sense, equivalent to \eqref{4.6}. The last two of these
give formulas occurring without proof in papers by Olsson.
We postpone a further corollary, a rewritten form of \eqref{4.6}
which is close to the Yoshida-Kita integral,
to Section \ref{s5}.

\begin{corollary}
\label{cor4.5}
For $\Re a,\Re(a-c_{2}+1),\Re(b_2+c_1-a)>0$
and $(x,y)\in\CC^2$ such that
$x\notin[1,\iy)$ and $y\notin (-\iy,0]$, we have
\begin{multline}
F_P(a,b_{1},b_{2},c_{1},c_{2};x,y)=
\frac{\Ga(c_1)\Ga(a+b_{2}-c_{2}+1)}
{\Ga(a)\Ga (a-c_{2}+1)\Ga(b_{2}+c_1-a)}\\
\times y^{-b_2}\int_0^1
u^{a-b_2-1} (1-u)^{b_2+c_1-a-1}(1-xu)^{-b_1}
\hyp21{b_2,b_2-c_2+1}{-a+b_2+c_1}{y^{-1}(1-u^{-1})}\,du.
\label{eq32}
\end{multline}
\end{corollary}

The integral representation \eqref{eq32} for $F_P$ is different
from two integral representations \cite[(43), (48)]{14} for $F_P$
in terms of ${}_2F_1$. By \eqref{eq33} and \eqref{2.2} we obtain
from \eqref{eq32} the specialization formula
\begin{equation}
F_P(a,b_1,b_2,a,c_2;0,y)=y^{-c_2+1}
\hyp21{b_2-c_2+1,a-c_2+1}{a+b_2-c_2+1}{1-y}\qquad(y\notin(-\iy,0]).
\label{eq34}
\end{equation}
\begin{corollary}
\label{cor4.2}
With the assumptions of Theorem \ref{th4.1} we have
\begin{multline}
F_P(a,b_{1},b_{2},c_{1},c_{2};x,y)=
\frac{\Ga(a+b_{2}-c_{2}+1)\Ga(c_{1})\Ga(b_{2}+c_{1}-a)}
{\Ga(b_{2}+c_{1})\Ga(a-c_{2}+1)\Ga( b_{2}-c_{2}+1)\Ga(c_{1}+c_{2}-a-1)}
\\
\times \int_{0}^{1}v^{-c_{2}}(1-v)^{c_{1}+c_{2}-a-2}
F_1\left(a,b_1,b_2,b_2+c_1;x,\frac{v-y}{v}\right) dv.
\label{4.7}
\end{multline}
\end{corollary}

Formula \eqref{4.7} was earlier stated without proof by
Olsson\footnote{In \cite[(65)]{14} the factor $(c_{1}+c_{2}-a+1)$
should be $\Ga(c_{1}+c_{2}-a+1)$.} \cite[(65)]{14}.

\begin{corollary}
\label{cor4.3}
With the assumptions of Theorem \ref{th4.1} and moreover $\Re y>\half$
we have
\begin{multline}
F_P(a,b_{1},b_{2},c_{1},c_{2};x,y)=
\frac{\Ga(c_{1})\Ga(a+b_{2}-c_{2}+1)}
{\Ga(a-c_{2}+1)\Ga(b_{2}+c_{1})}\,y^{-b_{2}}\\
\times
\sum_{j=0}^\iy\frac{(b_{2})_j(b_2+c_1-a)_j}{(b_{2}+c_{1})_j\,j!}
\hyp21{a,b_{1}}{b_{2}+c_{1}+j}x
\hyp21{-j,b_{2}-c_{2}+1}{b_2+c_1-a}{y^{-1}}.
\label{4.7a}
\end{multline}
\end{corollary}

Formula \eqref{4.7a}
was earlier stated without proof by Almstr\"{o}m \& Olsson
\cite[(22)]{18}. See also \cite[(55)]{14}.
\bLP
{\bf Proof of Theorem \ref{th4.1}.}
Denote the \RHS\ of \eqref{4.6} by
$I(a,b_1,b_2,c_1,c_2;x,y)$ and consider this under the given
constraints for $x,y$ and the parameters.
Then all occurring Gamma factors have argument with positive real part.
Because
\[
\const u\le \big|u+t(1-u)vy^{-1}\big|\le\max\big(|y|^{-1},1\big),
\]
uniformly for $y$ in compact subsets of $\CC$ outside $(-\iy,0]$,
we see that the double integral converges absolutely, and that
$I(x,y)$ is analytic for $(x,y)$ in the given region.
Now restrict $(x,y)$ to $-1<x<1$ and $y>0$. We will show that
$I(x,y)$ is then equal to the series \eqref{4.2}. Thus the functions
given by \eqref{4.2} and by \eqref{4.1} on their initial domains
will have analytic continuation to the domain of $I(x,y)$.

Expand the factor $(1-xu)^{-b_{1}}$
in the integrand as binomial series and
interchange the sum and the integrals by dominated convergence.
We obtain
\begin{equation}
I(a,b_{1},b_{2},c_{1},c_{2};x,y)=
\sum_{j=0}^\iy\frac{(b_{1})_j}{j!}\,x^j\,I_j(a,b_{2},c_{1},c_{2};y)
\label{4.9}
\end{equation}
with
\begin{multline*}
I_j(a,b_{2},c_{1},c_{2};y)=I_j(y)=
\frac{\Ga(a+b_{2}-c_{2}+1)\Ga(c_{1})}
{\Ga(a)\Ga (a-c_{2}+1)\Ga(b_{2}-c_{2}+1)\Ga(c_{1}+c_{2}-a-1)}\\
\times\int_0^1 v^{-c_{2}}(1-v)^{c_{1}+c_{2}-a-2}
\int_0^1 u^{a+j-1}(1-u)^{b_{2}+c_{1}-a-1}
\left(1-\frac{v-y}{v}u\right)^{-b_{2}}\,du\,dv.
\end{multline*}
The inner integral can be evaluated by \eqref{2.2}. Then
\begin{align*}
I_j(y)&=\frac{\Ga(a+b_{2}-c_{2}+1)\Ga(c_{1})\Ga(b_{2}+c_{1}-a)}
{\Ga(b_{2}+c_{1})\Ga(a-c_{2}+1)\Ga(b_{2}-c_{2}+1)\Ga(c_{1}+c_{2}-a-1)}\,
\frac{(a)_j}{(b_{2}+c_{1})_j}\\
&\qquad\qquad\qquad\qquad
\times\int_0^1 v^{-c_{2}}(1-v)^{c_{1}+c_{2}-a-2}
\hyp21{b_{2},a+j}{b_{2}+c_{1}+j}{\frac{v-y}{v}}dv\\
&=\frac{\Ga(a+b_{2}-c_{2}+1)\Ga(c_{1})\Ga(b_{2}+c_{1}-a)}
{\Ga(b_{2}+c_{1})\Ga(a-c_{2}+1)\Ga(b_{2}-c_{2}+1)\Ga(c_{1}+c_{2}-a-1)}
\frac{(a)_j}{(b_{2}+c_{1})_j}y^{1-c_2}\\
&\qquad\qquad\qquad\qquad
\times\int_{-\iy}^{1-y}(1-w)^{a-c_{1}}(1-y-w)^{c_{1}+c_{2}-a-2}
\hyp21{b_{2},a+j}{b_{2}+c_{1}+j}w\,dw.
\end{align*}
Now use \cite[(4.11)]{11}
\begin{multline*}
\int_{-\iy}^{x}(1-y)^{a+b-c}\hyp21{a,b}{c}y
\frac{(x-y)^{\mu -1}}{\Ga(\mu)}\,dy\\
=\frac{\Ga(c-a-\mu)}{\Ga(c-a)}\,
\frac{\Ga(c-b-\mu)}{\Ga(c-b)}\,
\frac{\Ga(c)}{\Ga(c-\mu)}\,
(1-x)^{a+b-c+\mu}\hyp21{a,b}{c-\mu}x,
\end{multline*}
with convergence conditions $x<1$, $\Re(c-a),\Re(c-b)>\Re(\mu)>0$,
in order to obtain
\[
I_j(y)=\frac{(a)_j(a-c_{2}+1)_j}{(c_1)_j(a+b_{2}-c_{2}+1)_j}
\hyp21{a+j,b_{2}}{a+b_{2}-c_{2}+1+j}{1-y},
\]
with $y>0$ and $\Re(c_1),\Re(a-c_2+1),\Re(c_1+c_2-a-1)>0$.
These constraints are implied by what we had already assumed.
Substitution in \eqref{4.9} and combination with \eqref{4.2} gives the
desired result \eqref{4.6}.\qed
\bLP
{\bf Proof of Corollary \ref{cor4.5}.}
First assume the conditions of Theorem \ref{th4.1}.
In the integrand of \eqref{4.6} extract the factor $u^{-b_2}$ from
$(u+(1-u)vy^{-1})^{-b_2}$ and then substitute Euler's integral
\eqref{2.1} for the $v$-integral.
This yields \eqref{eq32}.
The conditions on the parameters in Corollary \ref{cor4.5} are implied
by those in Theorem \ref{th4.1} and they are sufficient for
absolute convergence of the integral (use \cite[2.10(2)]{6}).
Hence, by analytic continuation in the parameters one can show 
that \eqref{eq32} holds under the given constraints.\qed
\bLP
{\bf Proof of Corollary \ref{cor4.2}.}
Rewriting of \eqref{4.6} gives
\begin{multline*}
F_P(a,b_{1},b_{2},c_{1},c_{2};x,y)=
\frac{\Ga(a+b_{2}-c_{2}+1)\Ga(c_{1})}
{\Ga(a)\Ga (a-c_{2}+1)\Ga(b_{2}-c_{2}+1)\Ga(c_{1}+c_{2}-a-1)}\\
\times\int_0^1 v^{-c_{2}}(1-v)^{c_{1}+c_{2}-a-2}
\int_0^1 u^{a-1}(1-u)^{b_{2}+c_{1}-a-1}(1-xu)^{-b_{1}}
\left(1-\frac{v-y}{v}u\right)^{-b_{2}}\,du\,dv.
\end{multline*}
The inner integral is an $F_1$ Appell function,
see \cite[5.8(5)]{6}.\qed
\bLP
{\bf Proof of Corollary \ref{cor4.3}.}
Assume moreover $\Re y>\half$. Then $|1-vy^{-1}|\le v$.
So we can expand
$(u+(1-u)vy^{-1})^{-b_{2}}=
\big(1-(1-y^{-1}v)(1-u)\big)^{-b_2}$ in \eqref{4.6}
as binomial series. By dominated convergence we obtain
\begin{multline*}
F_P(a,b_{1},b_{2},c_{1},c_{2};x,y)=
\frac{\Ga(a+b_{2}-c_{2}+1)\Ga(c_{1})}
{\Ga(a)\Ga (a-c_{2}+1)\Ga(b_{2}-c_{2}+1)\Ga(c_{1}+c_{2}-a-1)}\,
y^{-b_{2}}\\
\times \sum_{j=0}^\iy\frac{(b_{2})_j}{j!}
\int_0^1 u^{a-1}(1-u)^{b_{2}+c_{1}-a-1+j}(1-xu)^{-b_{1}}
\int_0^1 v^{b_{2}-c_{2}}(1-v)^{c_{1}+c_{2}-a-2}
(1-y^{-1}v)^j\,du\,dv
\end{multline*}
Then substitution of \eqref{2.2} or \eqref{2.4} gives the result.\qed
\section{The Yoshida-Kita integral for $H_2$}
\label{s5}
In this section we compare a rewritten version of the double integral
representation \eqref{4.6} for $F_P$ with Yoshida \cite[(0.10)]{17} and 
Kita \cite[p.57, item 9]{10}, who, on first superficial view,
seem to give an evaluation of this integral as an $H_2$ function.
The difference is caused because Kita's integral is defined quite
differently, by using homological methods involving twisted cycles,
while we are integrating classically over a domain in $\RR^2$.
We will consider a variant of our double integral \eqref{4.6}
in which the outer integral is taken over a
Pochhammer double loop. This turns out to have the same $H_2$
evaluation as Kita's integral. There is
a branch point inside the double loop caused by
an independent variable. This prevents shrinking of the double loop
outer integral such that one arrives at the double integral \eqref{4.6}.

For a better understanding
we first discuss a comparable case of a double loop integral version
of Euler's integral repesentation for ${}_2F_1$ which does not allow
immediate shrinking. It turns out that there
shrinking of the contour is possible
to a sum of two Euler type integrals, by which we arrive at a
standard three-term identity for~${}_2F_1$.

Then we do a similar,
but more complicated,
exercise for the double integral with
double loop outer integral which has an evaluation as an $H_2$ function.
We will keep the inner integral as an ordinary integral over $[0,1]$,
and we will evaluate it by Euler's integral \eqref{2.2}.
Then we are dealing with a double loop integral
having a Gauss hypergeometric function in the integrand.
After shrinking of the double loop this integral will split as a sum
of two ordinary integrals over intervals.
We will treat only one of the four occurring cases in more detail. 
It will turn out that both terms can be evaluated as an $F_P$ function
by means of \eqref{eq32}.
The resulting
three-term indentity will now be an identity given by Olsson \cite{14},
involving two $F_P$ terms and one $H_2$ term.

\subsection{Pochhammer's double loop integral
for the Gauss hypergeometic function}
\label{s5.1}
The notion of a {\em double loop} integral was introduced by Jordan
\cite[pp.~243--244]{22} (1887) and Pochhammer
\cite{24} (1890). This means an integral over a contour in
the complex $t$-plane starting at some point $t_0$
on the interval $(0,1)$, for which we assume that $\arg t_0=0$
and $\arg(1-t_0)=0$,
then making
successively a loop around 1 in the positive sense, around 0
in the positive sense, around 1 in the negative sense, and
around 0 in the negative sense, and finally returning
at $t_0$ (see picture, for instance, in
\cite[p.272]{6}).
The prototypical evaluation of such an integral was done by
Pochhammer \cite{24} for the generalized beta integral,
see also \cite[1.6(7)]{6}:
\begin{equation}
\frac1{(1-e^{2\pi ia})(1-e^{2\pi ib})}
\int^{(1+,0+,1-,0-)}t^{a-1}(1-t)^{b-1}\,dt=
\frac{\Ga(a)\Ga(b)}{\Ga(a+b)}\qquad(a,b\in\CC\backslash\ZZ).
\label{eq1}
\end{equation}
For $\Re a,\Re b>0$ it can be reduced to the classical beta integral
by shrinking the loops around 0 and 1.

Since the specifications for the loop are only of topological nature,
the actual choice of the loop allows much freedom. However,
a convenient standard form, in particular if we want to shrink the loop,
is to take $\ep\in(0,\half)$, define $C_\ep(z)$ as the circle
of radius $\ep$ around $z\in\CC$, and then let the loop be given
as follows. Start at the pont $\ep$ on the interval $[0,1]$.
First go to $1-\ep$ along the interval, then turn around 1
along $C_\ep(1)$ in the positive sense. Then go from $1-\ep$ to $\ep$
and turn around 0 along $C_\ep(0)$ in the positive sense.
Then go from $\ep$ to $1-\ep$ and turn around 1 along $C_\ep(1)$ 
in the negative sense. Finally go from $1-\ep$ to $\ep$
and turn around 0 along $C_\ep(0)$ in the negative sense.

As already observed by Pochhammer \cite{24}, not just the beta integral
can be extended to a wider parameter range by integrating over a
double loop, see \eqref{eq1}, but this similarly works for
Euler's integral representation \eqref{2.2}
for the Gauss hypergeometric function, see \cite[2.1(13)]{6}:
\begin{align}
\hyp21{a,b}cz&=\frac{\Ga(c)}{\Ga(b)\Ga(c-b)}\,
\frac1{(1-e^{2\pi ib})(1-e^{2\pi i(c-b)})}\nonu\\
&\qquad\qquad\times
\int^{(1+,0+,1-,0-)} t^{b-1}(1-t)^{c-b-1}(1-tz)^{-a}\,dt
\label{eq9}\\
&\quad(b,1-c,c-b\notin\ZZ_{>0},\;z^{-1}\notin[0,1]\;
\mbox{and $z^{-1}$ outside the double loop}).\nonu
\end{align}
The proof is by reduction to \eqref{eq1}, where we now assume that
$b$ and $c-b$ avoid the integers completely and that $|z|<1$.
Indeed, deform the double
loop such that $|tz|<1$ for $t$ on the loop (allowed because $z^{-1}$ is
outside the loop and $|z^{-1}|>1$) and expand $(1-tz)^{-a}$ as a
binomial series.
Note also that,
for $\Re b, \Re(c-b)>0$ and with $z^{-1}\notin[0,1]$ and fixed outside the
double loop, we can shrink the loops around 0 and 1
in order to arrive at \eqref{2.2}.

Now consider the double loop integral in \eqref{eq9}
under the assumption that $z\notin[0,\iy)$, that $z^{-1}$ is not
on the loop, and that there is a path connecting $z^{-1}$ with 0
which does not cross the double loop. Thus we consider
\begin{multline}
\int^{(1+,\{z^{-1},0\}+,1-,\{z^{-1},0\}-)} t^{b-1}(1-t)^{c-b-1}
(1-tz)^{-a}\,dt\\
=(-z)^{-a}\int^{(1+,\{z^{-1},0\}+,1-,\{z^{-1},0\}-)} t^{b-a-1}(1-t)^{c-b-1}
(1-t^{-1}z^{-1})^{-a}\,dt=:I(z).
\label{eq10}
\end{multline}
Assume that moreover $|z^{-1}|<1$ and that $b-a,c-b\notin\ZZ$.
Then we can deform the double
loop in $I(z)$ such that $|t^{-1}z^{-1}|<1$ for $t$ on the loop.
Binomial expansion, interchange of sum and integral by
dominated convergence, and application of \eqref{eq1} gives that
\begin{equation*}
I(z)=\frac{\Ga(b-a)\Ga(c-b)}{\Ga(c-a)}\,(1-e^{2\pi i(b-a)})(1-e^{2\pi i(c-b)})
(-z)^{-a}\sum_{k=0}^\iy\frac{(a)_k (a-c+1)_k}{(a-b+1)_k\,k!}\,z^{-k}.
\end{equation*}
Thus we have shown that
\begin{align}
&\frac{\Ga(c)\Ga(b-a)}{\Ga(b)\Ga(c-a)}\,(-z)^{-a}\,
\hyp21{a,a-c+1}{a-b+1}{z^{-1}}=
\frac{\Ga(c)}{\Ga(b)\Ga(c-b)}\,
\frac1{(1-e^{2\pi i(b-a)})(1-e^{2\pi i(c-b)})}\nonu\\
&\qquad\qquad\qquad\qquad\qquad
\times\int^{(1+,\{z^{-1},0\}+,1-,\{z^{-1},0\}-)}
t^{b-1}(1-t)^{c-b-1}(1-tz)^{-a}\,dt
\label{eq6}
\\
&\qquad
(b-a,c-b\notin\ZZ,\;c\notin\ZZ_{\le0},\;z\notin[0,\iy),\;
\mbox{$z^{-1}$ connected with 0 inside the double loop}).\nonu
\end{align}

\begin{figure}[!h]
\centering     
\includegraphics[width=10cm]{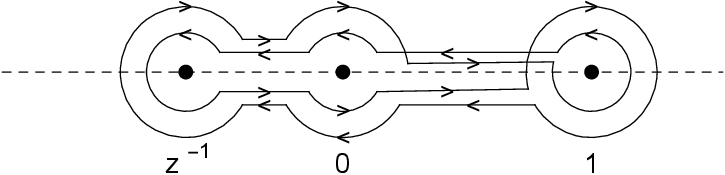} 
\caption{Double loop prepared for shrinking in the integral in \eqref{eq6}.}
\label{fig:3}
\end{figure}

It is also possible to shrink the double loop in \eqref{eq6} and thus
arrive at a sum of two classical Euler type integrals. For convenience
assume that $z<0$ and deform the double loop such that it consists
of segments $[z^{-1}+\ep,-\ep]$, $[\ep,1-\ep]$, circles of radius
$\ep$ around $z^{-1}$ and 1, and two half circles of radius $\ep$
around $0$, where each part is traversed various times and in
different directions (see Figure \ref{fig:3}).
In Table \ref{table:2} we list the arguments of the factors in the
integrand for $t$ on the various parts of the double loop.
\begin{table}[!h]
\centering
\begin{tabular}{l|c c c c c c c c}
&$\ri{(0,1)}$&$\li{(0,1)}$&$\li{(z^{-1},0)}$&$\ri{(z^{-1},0)}$&
$\ri{(0,1)}$&$\li{(0,1)}$&$\li{(z^{-1},0)}$&$\ri{(z^{-1},0)}$
\\[\smallskipamount]
\hline
$\arg(t)$&0&0&$\pi$&$\pi$&$2\pi$&$2\pi$&$\pi$&$\pi$\\
$\arg(1-t)$&0&$2\pi$&$2\pi$&$2\pi$&$2\pi$&0&0&0\\
$\arg(1-tz)$&0&0&0&$2\pi$&$2\pi$&$2\pi$&$2\pi$&0\\
\end{tabular}
\caption{Essential data while $t$ goes through the double loop.}
\label{table:2}
\end{table}

Also assume that $\Re a<1$ and $\Re c>\Re b>0$
with $b-a,c-b\notin\ZZ$, $c\notin\ZZ_{\le0}$.
Then, as $\ep\downarrow0$, the \RHS\ of \eqref{eq6} tends to
\begin{align}
&\frac{\Ga(c)}{\Ga(b)\Ga(c-b)}\left(
\int_0^1t^{b-1}(1-t)^{c-b-1}(1-tz)^{-a}\,dt\right.\nonu\\
&\qquad\qquad\qquad\qquad
\left.-\,\frac{\sin(\pi a)}{\sin(\pi(a-b))}
\int_{z^{-1}}^0 (-t)^{b-1}(1-t)^{c-b-1}(1-tz)^{-a}\,dt\right)\nonu\\
&=\hyp21{a,b}cz-\,\frac{\Ga(c)\Ga(a-b)\Ga(b-a+1)}
{\Ga(b)\Ga(c-b)\Ga(a)\Ga(1-a)}\,(-z)^{-b}
\int_0^1 s^{b-1}(1-s)^{-a} (1-z^{-1}s)^{c-b-1}\,ds\nonu\\
&=\hyp21{a,b}cz-\,\frac{\Ga(c)\Ga(a-b)}{\Ga(a)\Ga(c-b)}\,(-z)^{-b}
\hyp21{b,b-c+1}{b-a+1}{z^{-1}},\label{eq11}
\end{align}
where we used \eqref{2.2} and \cite[1.2(6)]{6}.
It follows that \eqref{eq11} is equal to \eqref{eq6}, by which we
have recovered \cite[2.10(2)]{6}.
\begin{remark}
\label{th5}
In the homological approach it is quite standard to prove
connection formulas as above. See for instance Mimachi
\cite{26} for the more general case dealing with ${}_3F_2$.
\end{remark}

\subsection{A classical analysis version of Kita's integral for $H_2$}
We start with another corollary to Theorem \ref{th4.1}.
\begin{corollary}
\label{cor4.4}
For $\Re a,\Re d,\Re(a+c),\Re(a+d),\Re(e-a-d)>0$
and $(x,y)\in\CC^2$ such that
$x\notin[1,\iy)$ and $y\notin[0,\iy)$, we have
\begin{align}
&\frac{\Ga(a+c)\Ga(a+d)}{\Ga(a+c+d)\Ga(a)}\,
(-y)^{-c} F_P(a+c,b,c,e,c-d+1;x,-y^{-1})
=\frac{\Ga(e)}{\Ga(e-a-d)\Ga(a)\Ga(d)}
\nonu\\
&\quad\times\int_0^1\int_0^1 u^{a-1} (1-u)^{e-a-1}
v^{d-1} (1-v)^{e-a-d-1} (1-ux)^{-b}
\big(1-(1-u)u^{-1}vy\big)^{-c}\,du\,dv
\label{4.8}\\
\noalign{\allowbreak}
&\qquad=\frac{\Ga(e)}{\Ga(e-a-d)\Ga(a)\Ga(d)}\,
\int_{u=0}^1 \int_{v=0}^{u^{-1}-1} u^{a+d-1} (1-ux)^{-b} v^{d-1}
(1-vy)^{-c}
\nonu\\
&\hskip9cm\times (1-u-uv)^{e-a-d-1}\,dv\,du
\label{4.8a}\\
&\qquad=\frac{\Ga(e)}{\Ga(e-a-d)\Ga(a)\Ga(d)}\,
\int_{u=0}^1 \int_{v=0}^{1-u} u^{a-1} (1-ux)^{-b} v^{d-1}
(1-u^{-1}vy)^{-c}
\nonu\\
&\hskip9cm\times (1-u-v)^{e-a-d-1}\,dv\,du.
\label{4.8b}
\end{align}
\end{corollary}
\noindent
{\bf Proof.}
Formula \eqref{4.8} is a trivial rewriting of \eqref{4.6}.
We went from \eqref{4.8} to \eqref{4.8a} by transforming
$(u,v)\mapsto\big(u,\frac{uv}{1-u}\big)$ and from \eqref{4.8a} to \eqref{4.8b}
by $(u,v)\mapsto (u,u^{-1}v)$.\qed
\bPP
Expression \eqref{4.8a} was evaluated
as $H_2(a,b,c,d,e;x,y)$ by Yoshida\footnote{In
\cite[(0.10)]{17} a power of $-(uv+u-1)$ should be taken,
because the integral is over a region with $uv+u-1\le0$.}
\cite[(0.10)]{17}. 
Kita \cite[p.57, item 9]{10} gave the same evaluation\footnote{In
\cite[p.57, item 9]{10} the second occurrence of $z_{1}$ should be $z_{2}$.
Also, on \cite[p.27, line 3]{10} the denominator $\Ga(c)$ should
be $\Ga(\gamma)$.} for \eqref{4.8b} as Yoshida gave for \eqref{4.8a},
but he specified
that integration should be understood by the homological approach,
which changes things drastically.
Dwork \& Loeser \cite[p.105]{5} also reproduce
Yoshida's formula, but without specification of the type of
integral.

We will now show that we can reproduce Kita's evaluation of \eqref{4.8b} as $H_2$
by working with classical integrals. For our purpose
it is easier to work with a variant of the
\RHS\ of \eqref{4.8}.

\begin{theorem}
Let $\Re d,\Re(e-a-d)>0$ and $a,e-a\notin\ZZ$. Then
\begin{multline}
H_2(a,b,c,d;e,x,y)=\frac{\Ga(e)}{\Ga(e-a-d)\Ga(a)\Ga(d)}\,
\frac1{(1-e^{2\pi ia})(1-e^{2\pi i(e-a)})}\int^{(1+,0+,1-,0-)}\\
\times u^{a-1} (1-u)^{e-a-1}(1-ux)^{-b}
\left(\int_0^1 v^{d-1} (1-v)^{e-a-d-1}
\big(1-(u^{-1}-1)vy\big)^{-c}\,dv\right)du.
\label{eq3}
\end{multline}
Here $(x,y)\in\CC^2$ and the double loop are such that
$ux\notin[1,\iy)$ and $(u^{-1}-1)y\notin[1,\iy)$ for $u$ on the double loop.
\end{theorem}
\Proof
First choose $x\in\CC$ with $|x|<1$. Then choose a double loop
$(1+,0+,1-,0-)$ (as defined in the beginning of Section \ref{s5.1})
such that $x^{-1}$ is not on the double loop and can be connected with
$\iy$ without crossing the double loop. Then
$u^{a-1} (1-u)^{e-a-1}(1-ux)^{-b}$
is well defined on the loop. Since the set
$\{(u^{-1}-1)v\mid u\in(1+,0+,1-,0-),\;v\in[0,1]\}$ is compact,
we will have $|(u^{-1}-1)vy|<1$ in the double integral \eqref{eq3} for
$|y|$ small enough. Thus, for such $x,y$, the whole integrand is
well defined, and so is the double integral.

Write the \RHS\ of \eqref{eq3} as
\[
\frac{\Ga(e)}{\Ga(e-a-d)\Ga(a)\Ga(d)}\,J.
\]
.In $J$ we can expand
$(1-ux)^{-b}$ and $\big(1-(1-u)u^{-1}vy\big)^{-c}$ as binomial series
and interchange the double sum with the double integral by
dominated convergence. We obtain that
\begin{align*}
J&=
\sum_{k,l=0}^\iy\frac{(b)_k(c)_l}{k!\,l!}\,x^k y^l\,
\int_0^1 v^{d+l-1}(1-v)^{e-a-d-1}\,dv
\int^{(1+,0+,1-,0-)}u^{a+k-l-1}(1-u)^{e-a+l-1}\,du\\
&=\sum_{k,l=0}^\iy\frac{(a)_{k-l} (b)_k (c)_l (d)_l}{(e)_k k!\,l!}\,
x^k y^l=H_2(a,b,c,d,e;x,y).
\end{align*}
Analytic continuation in $x$
and $y$ with simultaneous deformation of the double loop is possible
as long as $ux\notin[1,\iy)$ and
$(u^{-1}-1)y\notin[1,\iy)$ for $u$ on the double loop.\qed
\subsection{The integral \eqref{eq3} as a solution of the $F_2$ system
and analytic continuation}
We can rewrite \eqref{eq3} as
\begin{align}
&y^{-b_2} H_2(a-b_2,b_1,b_2,b_2-c_2+1,c_1;x,-y^{-1})
=\frac{\Ga(c_1)}{\Ga(-a+c_1+c_2-1)\Ga(a-b_2)\Ga(b_2-c_2+1)}\nonu
\\
&\times\frac{y^{-b_2}}{(1-e^{2\pi i(a-b_2)})(1-e^{2\pi i(-a+b_2+c_1)})}\,
\int^{(1+,0+,1-,0-)}
u^{a-b_2-1} (1-u)^{-a+b_2+c_1-1}(1-xu)^{-b_1}\nonu
\\
&\qquad\qquad\quad
\times\left(\int_0^1 v^{b_2-c_2} (1-v)^{-a+c_1+c_2-2}
\big(1-y^{-1}(1-u^{-1})v\big)^{-b_2}\,dv\right)du=: I(x,y),
\label{eq20}
\end{align}
where $\Re(b_2-c_1+1),\Re(-a+c_1+c_2-1)>0$ and
$a-b_2,-a+b_2+c_1\notin\ZZ$,
while $|x|<1$ and $|y|$ is large enough, with the same branch choice
for the two occurrences of $y^{-b_2}$.
By \cite[(14)]{14} the left-hand side of \eqref{eq20} is a
solution of the  system \cite[5.9(10)]{6} of pde's for $F_2$.
Actually, see also \eqref{3.1}, it is the unique solution
of this system of the form $y^{-b_2}\sum_{i,j=0}^\iy c_{i,j}x^i y^{-j}$
with $c_{0,0}=1$ which is analytic in a neigbourhood of $(x,y)=(0,\iy)$.

The {\em singular loci} of the $F_2$ system are the points around which
the local solution space has dimension less than the generic dimension four.
It is well known (and easily derived) that the singular loci are the lines
$x=0$, $x=1$, $y=0$, $y=1$, $x+y=1$ and (after one-point compactification
of the two $\CC$-factors of $\CC^2$) $x=\iy$, $y=\iy$.
We conclude that the
left-hand side of \eqref{eq20} can be certainly analytically extended
to a one-valued analytic function on the regions in $\RR^2$ given by
$y>1$, $1-y<x<1$ and by $x<1$, $y<0$. We will show that this can also be
read off from the \RHS\ $I(x,y)$ of \eqref{eq20}.

First we evaluate the inner integral in $I(x,y)$ as an Euler integral
\eqref{2.2}:
\begin{align}
&I(x,y)=\frac{\Ga(c_1)}{\Ga(a-b_2)\Ga(-a+b_2+c_1)}\,
\frac{y^{-b_2}}{(1-e^{2\pi i(a-b_2)})(1-e^{2\pi i(-a+b_2+c_1)})}\nonu\\
&\times
\int^{(1+,0+,1-,0-)}u^{a-b_2-1} (1-u)^{-a+b_2+c_1-1}(1-xu)^{-b_1}
\hyp21{b_2,b_2-c_2+1}{-a+b_2+c_1}{y^{-1}(1-u^{-1})}\,du.
\label{eq21}
\end{align}
Then the constraints
$\Re(b_2-c_1+1),\Re(-a+c_1+c_2-1)>0$ for \eqref{eq20} can be dropped,
while the other constraints are still kept.
For analytic continuation in \eqref{eq21} we can simultaneously bring
$(x,y)$ to a larger region and deform the double loop.
We will certainly have a one-valued analytic continuation as long as,
for $u$ in the double loop, $xu$ and $y^{-1}(1-u^{-1})$ do not cross
the cut $[1,\iy)$. For $(x,y)\in\RR^2$ and in the mentioned domains
$y>1$, $1-y<x<1$ or $x<1$, $y<0$ of unique analytic continuation for
$y^{-b_2}H_2(x,-y^{-1})$ we distinghuish four cases. We will see that
in each case the two cuts are no obstruction for the double loop.
\begin{enumerate}
\item
$x\le0$ and $x+y>1$. Then $u$ has to remain outside
$(-\iy,x^{-1}]$ and $[(1-y)^{-1},0]$, so $u$ can pass $\RR$ through the
holes $(1,\iy)$ and $(x^{-1},(1-y)^{-1})$ (left from 0).
\item
$0\le x<1$ and $y>1$. Then $u$ has to remain outside
$[(1-y)^{-1},0]$ and $[x^{-1},\iy)$, so $u$ can pass $\RR$ through the
holes $(1,x^{-1})$ and $(-\iy,(1-y)^{-1})$ (left from 0).
\item
$x\le0$ and $y<0$. Then $u$ has to remain outside
$(-\iy,x^{-1}]$ and $[0,(1-y)^{-1}]$, so $u$ can pass $\RR$ through the
holes $(1,\iy)$ and $(x^{-1},0)$.
\item
$0\le x<1$ and $y<0$. Then $u$ has to remain outside
$[0,(1-y)^{-1}]$ and $[x^{-1},\iy)$, so $u$ can pass $\RR$ through the
holes $(1,x^{-1})$ and $(-\iy,0)$.
\end{enumerate}
In a way depending on the case the double loop can be next shrunk,
analogous to our treatment of \eqref{eq6}.
As an illustration we will only treat the first case.

\subsection{Shrinking the double loop in \eqref{eq21}}
\begin{figure}[!ht]
\centering     
\includegraphics[width=8cm]{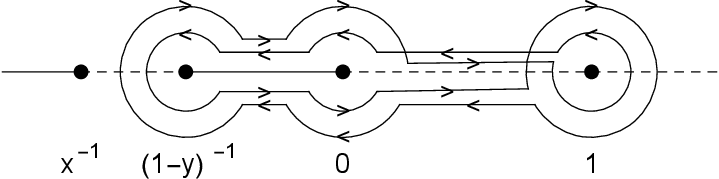} 
\caption{Double loop prepared for shrinking in the integral in \eqref{eq21}.}
\label{fig:4}
\end{figure}

\begin{figure}[!ht]
\centering
\begin{minipage}{6cm}   
\includegraphics[width=6cm]{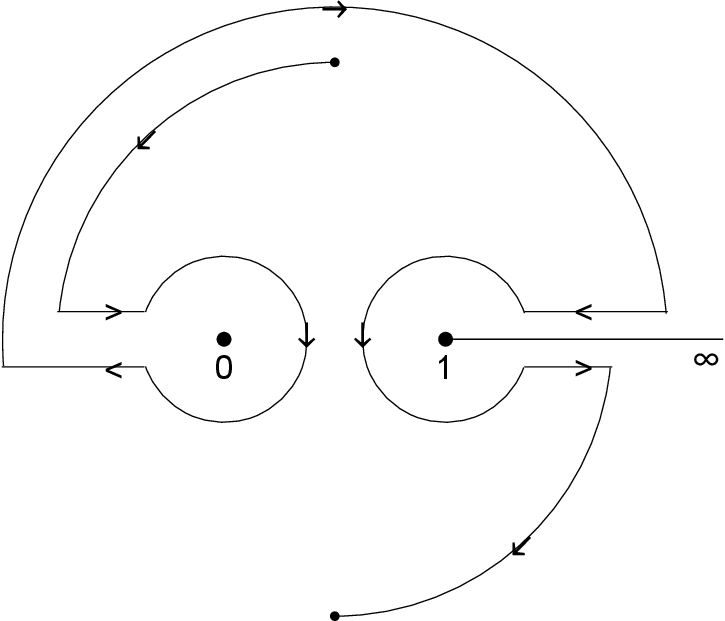} 
\end{minipage}
\quad
\begin{minipage}{6cm}   
\includegraphics[width=6cm]{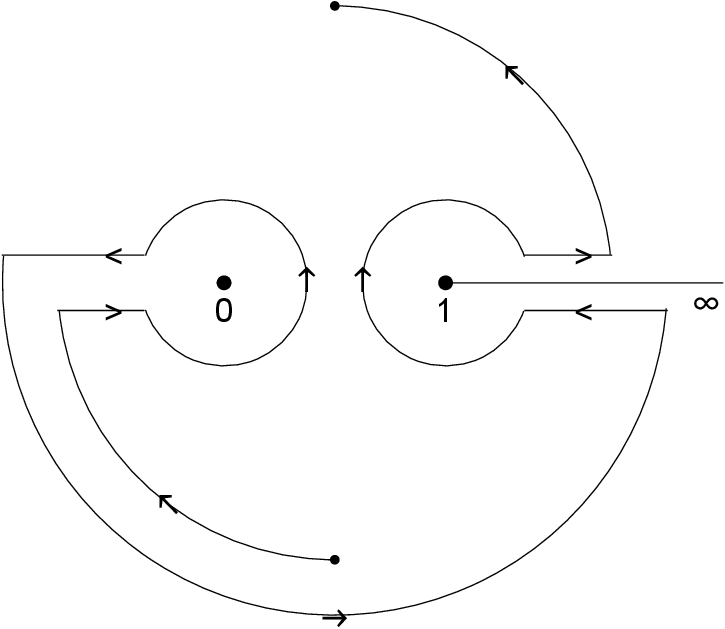} 
\end{minipage}4\caption{Orbits of $(1-u^{-1})y^{-1}$ during the first
and second part of the double loop.}
\label{fig:5}
\end{figure}

We will continue with case 1 above. Assume that $x<0$ and $x+y>1$.
Deform the double loop in \eqref{eq21} such that it consists
of segments $[(1-y)^{-1}+\ep,-\ep]$, $[\ep,1-\ep]$, circles of radius
$\ep$ around $(1-y)^{-1}$ and 1, and half circles of radius $\ep$
around 0, where each part is traversed various times and in
different directions (see Figure \ref{fig:4}).
One has to do good bookkeeping of the arguments of
$u$, $1-u$ and $y^{-1}(1-u^{-1})$ while $u$ goes through the double loop,
since these will affect the integrand.
See Figure \ref{fig:5} for the orbit of 
$(1-u^{-1})y^{-1}$ while $u$ goes through the double loop.
See Table \ref{table:1} for the essential data associated with the
double loop.
\begin{table}[!ht]
\centering
\begin{tabular}{r|c c c c}
&$\ri{(0,1)}$&$\li{(0,1)}$&$\li{((1-y)^{-1},0)+i0}$&
$\ri{((1-y)^{-1},0)-i0}$\\[\smallskipamount]
\hline
$\arg(u)=$&0&0&$\pi$&$\pi$\\
$\arg(1-u)=$&0&$2\pi$&$2\pi$&$2\pi$\\
$(1-u^{-1})y^{-1}\in$&$(-\iy,0)$&$(-\iy,0)$&$(1,\iy)+i0$&$(1,\iy)-i0$\\
\multicolumn{1}{r}{}&&&\\
&$\ri{(0,1)}$&$\li{(0,1)}$&
$\li{((1-y)^{-1},0)-i0}$&$\ri{((1-y)^{-1},0)+i0}$\\[\smallskipamount]
\hline
$\arg(u)=$&$2\pi$&$2\pi$&$\pi$&$\pi$\\
$\arg(1-u)=$&$2\pi$&0&0&0\\
$(1-u^{-1})y^{-1}\in$&$(-\iy,0)$&$(-\iy,0)$&$(1,\iy)-i0$&$(1,\iy)+i0$\\
\end{tabular}
\caption{Essential data while $u$ goes through the double loop.}
\label{table:1}
\end{table}
\newpage
We already had the constraints
$a-b_2,-a+b_2+c_1\notin\ZZ$. In order to allow the shrinking as
$\ep\downarrow0$ we also need the constraints
\begin{itemize}
\item
at $u=1$: $\Re(b_2+c_1-a)>0$;
\item
at $u=0$: $\Re a,\Re(a-c_2+1)>0$ (use \cite[2.10(2)]{6});
\item
at $u=(1-y)^{-1}$: $\Re(-a-b_1+c_1+c_2)>0$ (use \cite[2.10(1)]{6}).
\end{itemize}
In the limit for $\ep\downarrow0$ there will result the sum of an integral
over $u\in[0,1]$ and an integral over $u\in[(1-y)^{-1},0]$. For both integrals
we collect four instances of the integrand in \eqref{eq21} with branch
choices following Table \ref{table:1}. These four instances are added,
where a term has positive sign if the arrow in Table \ref{table:1} is in forward
direction, and negative sign for an arrow in negative direction.
We obtain
(also use \cite[1.2(6)]{6}):
\begin{equation}
I(x,y)=I_1(x,y)+I_2(x,y)
\label{eq25}
\end{equation}
with
\begin{multline}
I_1(x,y)=\frac{\Ga(c_1)}{\Ga(a-b_2)\Ga(-a+b_2+c_1)}\,
y^{-b_2}
\int_0^1
u^{a-b_2-1} (1-u)^{-a+b_2+c_1-1}(1-xu)^{-b_1}\\
\times\hyp21{b_2,b_2-c_2+1}{-a+b_2+c_1}{y^{-1}(1-u^{-1})}\,du
\label{eq22}
\end{multline}
and
\begin{multline}
I_2(x,y)=\frac{\Ga(c_1)\Ga(-a+b_2+1)}{\Ga(-a+b_2+c_1)}\,
\frac{y^{-b_2}}{2\pi i}
\int_{(1-y)^{-1}}^0
(-u)^{a-b_2-1} (1-u)^{-a+b_2+c_1-1}(1-xu)^{-b_1}\\
\times\left(\hyp21{b_2,b_2-c_2+1}{-a+b_2+c_1}{y^{-1}(1-u^{-1})+i0}
-\hyp21{b_2,b_2-c_2+1}{-a+b_2+c_1}{y^{-1}(1-u^{-1})-i0}\right)du.
\label{eq23}
\end{multline}

In \eqref{eq22} we can relax the parameter constraints to those
given above at $u=0$ and $u=1$.
By \eqref{eq32} we see that
\begin{equation}
I_1(x,y)=
\frac{\Ga(a)\Ga(a-c_2+1)}{\Ga(a+b_2-c_2+1)\Ga(a-b_2)}\,
F_P(a,b_{1},b_{2},c_{1},c_{2};x,y),
\label{eq30}
\end{equation}

In \eqref{eq23} we can write the difference of the two ${}_2F_1$
as one ${}_2F_1$ in view of \cite[2.10(1)]{6}, and next apply
Euler's transformation \eqref{2.6}:
\begin{align*}
&\frac1{2\pi i}\left(\hyp21{a,b}c{x+i0}-\hyp21{a,b}c{x-i0}\right)\\
&\qquad=\frac{\Ga(c)}{\Ga(a)\Ga(b)\Ga(c-a-b+1)}\,(x-1)^{c-a-b}\,
\hyp21{c-a,c-b}{c-a-b+1}{1-x}\\
&\qquad=\frac{\Ga(c)}{\Ga(a)\Ga(b)\Ga(c-a-b+1)}\,x^{1-c}(x-1)^{c-a-b}\,
\hyp21{1-a,1-b}{c-a-b+1}{1-x}\quad(x\in(1,\iy)).
\end{align*}
Then \eqref{eq23} takes the form
\begin{align}
&I_2(x,y)=
\frac{\Ga(c_1)\Ga(b_2-a+1)}{\Ga(b_2)\Ga(b_2-c_2+1)\Ga(c_1+c_2-a-b_2)}\,
y^{b_2-c_2}\int_{(1-y)^{-1}}^0(-u)^{a+b_2-c_2-1}\nonu\\
&\quad\times
\big(1+u(y-1)\big)^{c_1+c_2-a-b_2-1}(1-xu)^{-b_1}
\hyp21{c_2-b_2,1-b_2}{c_1+c_2-a-b_2}{1-y^{-1}(1-u^{-1})}\,du\nonu\\
\noalign{\allowbreak}
&=\frac{\Ga(c_1)\Ga(b_2-a+1)}{\Ga(b_2)\Ga(b_2-c_2+1)\Ga(c_1+c_2-a-b_2)}
\,(y-1)^{-a}\left(\frac y{y-1}\right)^{b_2-c_2}
\int_0^1 w^{a+b_2-c_2-1}\nonu\\
&\quad\times(1-w)^{c_1+c_2-a-b_2-1}\left(1-\frac{xw}{1-y}\right)^{-b_1}
\hyp21{c_2-b_2,1-b_2}{c_1+c_2-a-b_2}{\frac{y-1}y(1-w^{-1})}\,dw\nonu\\
&=\frac{\Ga(b_2-a+1)\Ga(a)\Ga(a-c_2+1)}
{\Ga(b_2)\Ga(a-b_2+1)\Ga(b_2-c_2+1)}\,(y-1)^{-a}\,
F_P\left(a,b_1,c_2-b_2,c_1,c_2;\frac x{1-y},\frac y{y-1}\right).
\label{eq24}
\end{align}
Here we have made the change of integration variable $u=w(1-y)^{-1}$
in the second equality, and we have used \eqref{eq32} in the
third equality. Furthermore, we can relax
in \eqref{eq24} the parameter constraints to those
given above at $u=0$ and $u=(1-y)^{-1}$. These are in agreement with
the parameter constraints in \eqref{eq32} (with $b_2$ replaced
by $c_2-b_2$).
\begin{remark}
The identity \eqref{eq25} with \eqref{eq20}, \eqref{eq22},
\eqref{eq24} substituted agrees with Olsson's identity \cite[(53)]{14}.
\end{remark}
\begin{remark}
If in the present section from \eqref{eq21} until here
we make the specializations $x=0$, $c_1=a$ then
the ${}_2F_1$ in \eqref{eq21} specializes to
$\big(1-y^{-1}(1-u^{-1})\big)^{c_2-b_2-1}$ by \eqref{eq33}.
The subsequent decomposition \eqref{eq25}
(i.e., $I(0,y)=I_1(0,y)+I_2(0,y)$) gives then another example of what we
wrote in
Section \ref{s5.1} starting with  \eqref{eq10}. There results
a formula connecting three Gauss hypergometric functions
(also seen by applying \eqref{eq35} and \eqref{eq34} to \eqref{eq20}
and \eqref{eq30}, \eqref{eq24}, respectively) which is essentially the
connection formula \cite[2.10(3)]{6}.
\end{remark}
\begin{remark}
Kato \cite{27} also discusses the solutions of the $F_2$ system at the
various singular points, but he does not refer to Olsson \cite{14},
and he gives, apart from $F_2$, few explicit solutions. However,
the first double power series in his Remark on p.330 can be
recognized as $F_P(a,b_2,b_1,c_2,c_1;y,x)$ expanded as \eqref{4.4}.
The three solutions $I(x,y)$, $I_1(x,y)$ and $I_2(x,y)$ 
of the $F_2$ system occurring above in \eqref{eq20}, \eqref{eq30},
\eqref{eq24}, respectively, can be matched up to a constant factor
with the following solutions in \cite{27}:
$I_1(x,y)$ with the first solution in \S3.2.1,
$I_2(x,y)$ with the first solution in \S3.3.5,
and $I(x,y)$ with the third solution in \S3.3.5.
Kato gives many connection formulas, but apparently not the
one connecting $I$, $I_1$ and $I_2$.
\end{remark}
\begin{remark}
Kita's integral representations
for Appell $F_2$ and $F_3$ \cite[p.56, items 2.2 and 3]{10},
which go back to Hattori \& Kimura \cite{25}, and where the double
integral is defined by the homological approach, remain valid
if the double integral is defined classically, but then further
constraints on the parameters are needed.
It would be interesting to consider also the other new double integral
representations given by Yoshida \cite[(0.6)--(0.9)]{17},
which come back in Kita \cite{10}
(for Yoshida's alternative $F_2$-integral \cite[(0.6)]{17} see also
Yoshida \cite[p.71]{21}).
Are these easy cases like the double integrals
for $F_2$ and $F_3$ or would their classical counterparts be a
sum of two integral representations?
\end{remark}
\section{Double integrals for solutions of the $F_2$ system}
\label{s6}
Olsson's list \cite[p.1289, Table I]{14} of
solutions of the system of PDEs
\cite[5.9(10)]{6} for $F_2$ includes, besides $F_2$ itself,
for instance functions expressible in terms of $F_3$, of $H_2$, and of $F_P$.
It is interesting to rewrite the well-known double
integral representations\footnote{In
\cite[5.8(3)]{6}, in the exponent of $(1-u-v)$, one should replace
$-\gamma$ by $\gamma$.}
\cite[5.8(2),(3)]{6} for $F_2$ and $F_3$
and our double integral representations \eqref{3.5} for $H_2$ and
\eqref{4.6} for $F_P$ as integral representations for solutions
of the $F_2$ system. We obtain:
\begin{align}
&F_2(a,b_1,b_2,c_1,c_2;x,y)=
\frac{\Ga(c_1)\Ga(c_2)}{\Ga(b_1)\Ga(b_2)\Ga(c_1-b_1)\Ga(c_2-b_2)}
\nonu\\
&\qquad\times x^{1-c_1}y^{1-c_2}
\int_{u=0}^x\int_{v=0}^y u^{b_1-1} v^{b_2-1} (1-u-v)^{-a}
(x-u)^{c_1-b_1-1}(y-v)^{c_2-b_2-1}\,dv\,du
\nonu\\
&\qquad\qquad\qquad\qquad\qquad
(0<x<1-y<1,\;\Re b_1, \Re b_2, \Re(c_1-b_1), \Re(c_2-b_2)>0),
\label{5.1}
\mLP
\noalign{\allowbreak}
&x^{-b_1} y^{-b_2}
F_3(b_1,b_2,1+b_1-c_1,1+b_2-c_2,b_1+b_2-a+1;x^{-1},y^{-1})
=\frac{\Ga(b_1+b_2-a+1)}{\Ga(b_1)\Ga(b_2)\Ga(1-a)}
\nonu\\
&\qquad\qquad\times x^{1-c_1}y^{1-c_2}
\int_{v=0}^1\int_{u=0}^{1-v} u^{b_1-1} v^{b_2-1} (1-u-v)^{-a}
(x-u)^{c_1-b_1-1}(y-v)^{c_2-b_2-1}\,du\,dv
\nonu\\
&\qquad\qquad\qquad\qquad\qquad\qquad\qquad\qquad\qquad\qquad
(x,y>1,\;\Re b_1,\Re b_2,\Re(1-a)>0),
\label{5.2}
\mLP
\noalign{\allowbreak}
&y^{-b_2} H_2(a-b_2,b_1,b_2-c_2+1,b_2,c_1;x,-y^{-1})
=\frac{\Ga(b_2-a+1)\Ga(c_1)}{\Ga(1-a)\Ga(b_2)\Ga(b_1)\Ga(c_1-b_1)}\,
x^{1-c_1}y^{1-c_2}
\nonu\\
&\qquad\qquad\qquad\times
\int_{u=0}^x\int_{v=0}^{1-u} u^{b_1-1} v^{b_2-1} 
(1-u-v)^{-a} (x-u)^{c_1-b_1-1} (y-v)^{c_2-b_2-1}\,dv\,du
\nonu\\
&\qquad\qquad\qquad\qquad\qquad
(0<x<1,\;y>1,\;\Re(1-a),\Re b_2, \Re b_1, \Re(c_1-b_1)>0),
\label{5.3}
\mLP
\noalign{\allowbreak}
&|x|^{1-c_1} y^{1-c_2} F_P(a-c_1-c_2+2,b_1-c_1+1,b_2-c_2+1,2-c_1,2-c_2;x,y)
\nonu\\
&\quad=\frac{\Ga(a+b_2-c_1-c_2+2)\Ga(2-c_1)}
{\Ga(a-c_1-c_2+2)\Ga(a-c_1+1)\Ga(b_2)\Ga(1-a)}\,
|x|^{1-c_1}y^{1-c_2}
\int_{u=1}^\iy
\int_{v=1-u}^0 u^{b_1-1}(-v)^{b_2-1}
\nonu\\
&\qquad\times (u+v-1)^{-a}(u-x)^{c_1-b_1-1}
(y-v)^{c_2-b_2-1}\,dv\,du\quad
\big(x\in(-\iy,0)\cup(0,1),y>0;
\nonu\\
&\qquad\quad\Re a<1,\;\Re b_2>0,\;
\Re(a-c_1-c_2+2),\Re(a-c_1+1),\Re b_2,\Re(1-a)>0\big).
\label{5.4}
\end{align}
Here we have restricted to domains of $(x,y)$ in $\RR^2$.
For going from \eqref{4.6} to \eqref{5.4} make the change of
integration variables
$(u,v)\to(u^{-1},v/(1-u))$.

Note that the \RHS s of \eqref{5.1}--\eqref{5.4} all have
the form
\begin{equation*}
\const|x|^{1-c_1}|y|^{1-c_2}\int\!\!\int 
|u|^{b_1-1} |v|^{b_2-1} |1-u-v|^{-a}
|x-u|^{c_1-b_1-1}|y-v|^{c_2-b_2-1}\,dv\,du,
\end{equation*}
where the integral is over a suitable $(u,v)$ region.
This will enable uniform proofs of such identities by showing that the
\RHS s are solutions of the system of pde's for $F_2$,
quite similar as was done in \cite{11} in the one-variable case. 
In Chapter 6 of the PhD Thesis by the first
author\footnote{E. Diekema, The fractional orthogonal derivative for
functions of one and two
variables, PhD Thesis, University of Amsterdam, 2018;
\url{http://hdl.handle.net/11245.1/a6ed8a3f-0831-4f9a-9476-2661ec7e1f92}}
the double integrals
\eqref{5.1}--\eqref{5.3} appear
as integral kernels in certain approximations
of two-dimensional fractional integrals.

\paragraph{Acknowledgements}
 The authors are very grateful to the referees. Their thoughtful
comments about
an earlier version of Section 5 (The Yoshida-Kita integral for $H_2$)
helped us greatly to do a complete rewriting and extension of
this section, providing a much deeper insight about the relationship
between
the Yoshida-Kita integral for $H_2$ on the one hand and our double
integral for $F_P$ on the other.
Comments about the subsequent version indicated an annoying error
and led to
a shortening of the proof of \eqref{eq24}. Many smaller details
could also be improved and important references could be added
as a result of these comments.
\vskip0.5cm
\begin{footnotesize}
\noindent
E. Diekema, Kooikersdreef 620, 7328 BS Apeldoorn, The Netherlands;
	\sPP
	email: {\tt e.diekema@gmail.com}
	\bLP
	T. H. Koornwinder, Korteweg-de Vries Institute, University of
	Amsterdam,\sPP
	P.O.\ Box 94248, 1090 GE Amsterdam, The Netherlands;
	\sPP
	email: {\tt T.H.Koornwinder@uva.nl}
\end{footnotesize}
	
\end{document}